\documentclass{amsart}

\newtheorem{theorem}{Theorem}[section]
\newtheorem{lemma}[theorem]{Lemma}
\newtheorem{conj}[theorem]{Conjecture}
\newtheorem{cor}[theorem]{Corollary}
\newtheorem{quest}[theorem]{Question}
\theoremstyle{definition}

\theoremstyle{remark}

\numberwithin{equation}{section}

\begin{document}

\title{A note on simple domains of GK dimension two}

\author{Jason P. Bell}
\address{Department of Mathematics, Simon Fraser University, 
8888 University Drive, Burnaby, BC, Canada,
V5A 1S6}
\email{jpb@math.sfu.ca}

\thanks{The author thanks NSERC for its generous support.}

\date{January 1, 2007.}


\begin{abstract}
Let $k$ be a field.  We show that a finitely generated simple Goldie $k$-algebra of quadratic growth is noetherian and has Krull dimension $1$.  In particular, if $A$ is a finitely generated simple domain of quadratic growth then $A$ is noetherian and by a result of Stafford every right and left ideal is generated by at most two elements.  We conclude by posing questions and giving examples in which we consider what happens when the hypotheses are relaxed.
\end{abstract}

\maketitle
\bibliographystyle{plain}

\section{Introduction}
We study finitely generated algebras of quadratic growth.  Given a field $k$ and a finitely generated $k$-algebra $A$, a $k$-subspace $V$ of $A$ is called a \emph{frame} of $A$ if $V$ is finite dimensional, $1\in V$, and $V$ generates $A$ as a $k$-algebra.
We say that $A$ has \emph{quadratic growth} if there exist a frame $V$ of $A$ and constants $C_1,C_2>0$ such that
$$C_1 n^2 \ \le \ {\rm dim}_k(V^n) \ \le \ C_2n^2\qquad {\rm for~all}~n\ge 1.$$
We note that an algebra of quadratic growth has GK dimension $2$.  More generally, the GK dimension of a finitely generated $k$-algebra $A$ is defined to be
$${\rm GKdim}(A) \ = \ \limsup_{n\rightarrow\infty} \log\bigg({\rm dim}(V^n)\bigg)\bigg/\log\, n,$$
where $V$ is a frame of $A$.  While algebras of quadratic growth have GK dimension $2$, it is not the case that an algebra of GK dimension $2$ necessarily has quadratic growth.  

GK dimension can be viewed as a noncommutative analogue of Krull dimension in the following sense: if $A$ is a finitely generated commutative $k$-algebra then the Krull dimension of $A$ and the GK dimension of $A$ coincide.  Krull dimension can also be extended to the noncommutative setting in a natural way by considering the so-called \emph{deviation} of the poset of left ideals of a ring (see McConnell and Robson \cite[Chapter 6]{MR} for the actual definition of Krull dimension).   It is in this more general sense that we use the term Krull dimension in this paper.  For more information on GK dimension we refer the reader to Krause and Lenagan \cite{KL}.

We study simple Goldie algebras of quadratic growth.  A ring $R$ is \emph{Goldie} if it satisfies the following two conditions:
\begin{enumerate}
\item $R$ does not contain an infinite direct sum of nonzero right ideals;
\item $R$ satisfies the ascending chain condition on right annihilators.
\end{enumerate}
Technically, this property is refered to as being \emph{right Goldie} (left Goldie is defined analogously), but a ring is right Goldie if and only if it is left Goldie.  

The Goldie condition is named after Alfred Goldie, who defined this property in his study of noncommutative localization.  Just as a commutative domain has a field of fractions, a prime Goldie ring $R$ has an Artinian ring of quotients $Q(R)$.  We note that noetherian algebras are always Goldie, but there exists non-noetherian Goldie algebras, for example a polynomial ring in infinitely many variables over a field is Goldie but not noetherian.  

Our main result is the following theorem, which shows that the Goldie and noetherian properties are equivalent for simple algebras of quadratic growth.  
\begin{theorem} \label{thm: main} Let $k$ be a field and let $A$ be a finitely generated simple Goldie $k$-algebra of quadratic growth.  Then $A$ is a noetherian algebra of Krull dimension $1$.  \end{theorem}
A result of Jategaonkar (cf. Krause and Lenagan \cite[Proposition 4.13]{KL}) shows that a domain of finite GK dimension is Goldie.  Hence we get the following corollary.
\begin{cor}  \label{cor: main} Let $k$ be a field and let $A$ be a finitely generated simple $k$-algebra of quadratic growth that is a domain.  Then $A$ is a noetherian algebra of Krull dimension $1$.  In particular, every left and right ideal can be generated by two elements.
\end{cor}
The fact that every right and left ideal of a simple noetherian algebra of Krull dimension one is generated by at most two elements follows from a result of Stafford \cite{St}.

The outline of the paper is as follows.  In section 2, we give proofs of Theorem \ref{thm: main} and Corollary \ref{cor: main}.  In section 3, we give some examples in which we consider what happens when various hypotheses are relaxed and give some concluding remarks.

\section{Proofs}
In this section we give the proofs of Theorem \ref{thm: main} and Corollary \ref{cor: main}.  We begin with a few estimates for algebras of quadratic growth. 
\begin{lemma} Let $K$ be a field, let $A$ be a finitely generated $K$-algebra with quadratic growth, and let $m$ be a positive integer.  If $V$ is a finite dimensional $K$-subspace of $A$ that generates $A$ as a $K$-algebra then there is a positive constant $C=C(m,V)$ such that
$${\rm dim}(V^n) - {\rm dim}(V^{n-m}) \ \le \ Cn$$ for infinitely many $n\ge 1$.
\end{lemma}
\begin{proof} By assumption, $A$ has quadratic growth and hence there is a positive constant $C_0$ such that
\begin{equation} {\rm dim}(V^n) \ \le \ C_0 n^2 \qquad {\rm for~}n ~{\rm sufficiently~large}.\end{equation}
We take \begin{equation}
C = 4 C_0 m.
\end{equation}
Suppose that there is a number $N$ such that $${\rm dim}(V^n) - {\rm dim}(V^{n-m}) > Cn$$ for all $n>N$.  Then 
\begin{eqnarray*}
{\rm dim}(V^{nm+N}) &\ge & {\rm dim}(V^{nm+N})-{\rm dim}(V^{N}) \\
& = & \sum_{j=0}^{n-1} {\rm dim}(V^{N+jm+m})-{\rm dim}(V^{N+jm})\\
& \ge & \sum_{j=0}^{n-1} C (N+jm+m) \\
&=& C Nn + Cm (n+1)n/2 \\
&\ge & 2 C_0 n^2m^2 \\
& > & C_0 (nm+N)^2 \qquad ({\rm for~}n~{\rm sufficiently~ large}).
\end{eqnarray*}
This is a contradiction.  The result follows. \end{proof}
Our next lemma gives an estimate for the growth of $A/I$ when $I$ contains a regular element $a$.  We recall that an element $a$ is \emph{regular} if its right and left annihilators are both $(0)$.
  \begin{lemma} Let $k$ be a field and let $A$ be a prime finitely generated $k$-algebra with quadratic growth. Let $V$ be a finite dimensional generating subspace of $A$ that contains $1$.  If $I$ is a right ideal of $A$ and $I$ contains a regular element $a$ then there is a positive constant $C$ such that
  $${\rm dim}(V^n/V^n\cap I) \ \le \ Cn$$ for infinitely many $n$. \label{lemma: 2}
  \end{lemma}
  \begin{proof}  By assumption $A$ has quadratic growth and hence there is a positive constant $C_0$ such that
$${\rm dim}(V^n) \ \le \ C_0n^2 \qquad {\rm for~}n\ge 1.$$
Let $I$ be a nonzero left ideal of $A$ and suppose there is a regular element $a\in I$.  Then $a\in V^m$ for some natural number $m$.  Then
\begin{eqnarray*}
{\rm dim}\left(V^n/(V^n\cap I)\right) &=& {\rm dim}(V^n) - {\rm dim}(V^n\cap I) \\
&\le & {\rm dim}(V^n) - {\rm dim}(V^n \cap aA)\\
&\le & {\rm dim}(V^n) - {\rm dim}(aV^{n-m})\\
&=& {\rm dim}(V^n) - {\rm dim}(V^{n-m}).\end{eqnarray*}
By the preceding lemma there is a positive constant $C$ such that
$${\rm dim}(V^n) - {\rm dim}(V^{n-m}) \ \le \ Cn$$ for infinitely many $n$.  The result follows. \end{proof}
\begin{lemma} Let $K$ be a field and let $A$ be an infinite simple affine $K$-algebra. Let $V$ be a finite dimensional generating subspace of $A$ which contains $1$.  If $I\subseteq J$ are two right ideals in $A$ and $I\not = J$ then there is $m>0$ such that $${\rm dim}(V^n\cap J/V^n\cap I) \ge n-m$$ for all $n\ge 0$. \label{lemma: 3}
\end{lemma}
\begin{proof} Pick $a\in J\setminus I$.   Without loss of generality, we may assume $J=I+aA$.  
Note that if $I+aV^m \supseteq I+aV^{m+1},$ then $I+aA\subseteq I+aV^m$, and hence $J/I$ is finite dimensional as a $K$-vector space.  But this says that $A$ has a finite dimensional nonzero right $A$-module, which is impossible since $A$ is infinite and simple.  It follows that
$${\rm dim}(I+aV^n/I) \ge n$$ for all $n\ge 0$.   Since $V$ is a generating subspace, there is some $m$ such that $a\in V^m$.  Hence
$${\rm dim}\left(J\cap V^{n+m}/I\cap V^{n+m}\right) \ \ge \ {\rm dim}\left((I+aV^{n})/I\right) \ \ge \ n.$$
The result follows. \end{proof}
 
{\em Proof of Theorem \ref{thm: main}.}
Let $V$ be a finite dimensional $K$-subspace of $A$ that generates $A$ as a $K$-algebra.  By assumption $A$ has quadratic growth and hence there is a positive constant $C$ such that
$${\rm dim}(V^n) \ \le \ Cn^2 \qquad {\rm for~}n\ge 1.$$
Let $I$ be a nonzero right ideal.  
By Lemma \ref{lemma: 2} there is a positive constant $C=C(I,V)$ such that
$${\rm dim}\left(V^n/(V^n\cap I)\right) \ \le \ Cn$$ for infinitely many $n$.  Let $m$ denote the \emph{uniform dimension} of $A$ (see McConnell and Robson \cite[2.2.10]{MR}).
We claim that if $$I=I_0\subseteq I_1 \subseteq I_2 \subseteq \cdots \subseteq I_e$$ is an ascending chain of ideals in which each containment is proper then $e\le mC$.  If not, there must exist $i$ and $j$ such that $d:=j-i>C$ and the uniform dimension of $A/I_i$ is the same as the uniform dimension of $A/I_{j}$.  Thus there exists a right ideal $J$ such that $I_{\ell}\cap J=(0)$ for $i\le \ell \le j$ and $I_{\ell}+J$ is essential (that is, it contains a regular element) for $i\le \ell \le j$.  Let
$$J_{\ell} = I_{i+\ell}+J$$ for $0\le \ell \le d$.

By Lemma \ref{lemma: 3} there are positive integers $m_1, m_2,\ldots m_d$ such that
$${\rm dim}(V^n\cap J_j/V^n\cap J_{j-1}) \ \ge \ n-m_j$$ for all $n$.
Hence
\begin{eqnarray*}
{\rm dim}(V^n) &\ge & {\rm dim}(V^n \cap J_d) \\
&=& {\rm dim}(V^n\cap J_0) + \sum_{i=1}^d {\rm dim}(V^n \cap J_i/V^n\cap J_{i-1}) \\
&\ge & {\rm dim}(V^n\cap J_0) +\sum_{i=1}^d (n-m_i)\\
&\ge & {\rm dim}(V^n\cap J_0) + dn - (m_1+\cdots +m_d).\end{eqnarray*}
Thus
$${\rm dim}(V^n) - {\rm dim}(V^n\cap J_0) \ \ge \ dn - (m_1+\cdots +m_d)$$ for all $n$.
But $J_0$ contains a regular element and hence for infinitely many $n$ we have
$${\rm dim}(V^n) - {\rm dim}(V^n\cap J_0) \ \le \ Cn$$ by Lemma \ref{lemma: 2}.
Thus
$d\le C$, since $m_1,\ldots ,m_d$ are fixed.   It follows that $A$ is right noetherian.  By looking at the opposite ring, we see $A$ is also left noetherian.  To see that $A$ has Krull dimension $1$, note that if 
$J\supset I$, and $I$ is nonzero, then any chain of submodules of $J/I$ in which containments are proper must have length at most $C=mC(I,V)$.  In particular, $J/I$ is artinian.  Hence $A$ has Krull dimension $1$.  By a result of Stafford, every right and left ideal of $A$ is generated by at most $2$ elements. \qed
\vskip 2mm
{\em Proof of Corollary \ref{cor: main}.} This follows immediately from Theorem \ref{thm: main} and Proposition 4.13 of Krause and Lenagan \cite{KL}. \qed \section{Concluding remarks}
We note that the finitely generated hypothesis is necessary. 
Let $$B=k[t,t^{1/2},t^{1/4},\ldots].$$  And let $A=B[x;\delta]$, where $\delta$ is differentiation with respect to $t$.  Then $A$ is a free $B$-module and since $B$ is not noetherian, neither is $A$.  We note that $A$ is a simple domain and any finitely generated subalgebra of $A$ has quadratic growth.

\[A = \left( \begin{array}{cc} k[x,y,y^{-1}] & k[x,y,y^{-1}]x \\
k[x,y,y^{-1}] & k[y]+k[x,y,y^{-1}]x \end{array} \right). \]
We note that $A$ is a finitely generated non-noetherian prime Goldie $k$-algebra (cf. McConnell and Robson \cite[13.10.2]{MR}).  It is straightforward to check that $A$ has quadratic growth.  
A non-PI example of a non-noetherian finitely generated simple Goldie algebra of quadratic growth can be constructed by using a field $k$ of characteristic $0$ and using $k_q[x,y,y^{-1}]$ in $A$ instead of  $k[x,y,y^{-1}]$, where $k_q[x,y,y^{-1}]$ is a localization of the \emph{quantum plane}; that is $xy=qyx$, with $q$ not a root of unity.  This type of construction cannot be made over a finite field.  Smoktunowicz \cite{Sm2} has shown that in many cases, a finitely generated prime algebra over a finite field satisfies a polynomial identity.  We thus make the following conjecture.
\begin{conj} Let $k$ be a finite field and let $A$ be a finitely generated $k$-algebra of quadratic growth.  Then $A$ is not simple.
\end{conj}

Artin and Stafford showed in the course of their description of graded domains of GK dimension two that non-PI graded domains of quadratic growth are noetherian.  It is an interesting question as to whether a similar description of the quotient division algebras of simple domains of quadratic growth can be given.  Artin and Stafford also made the conjecture that a ``gap'' should exist in the possible GK dimensions a graded domain can attain.  Smoktunowicz \cite{Sm3} proved this conjecture.  It would be interesting if a similar gap could be attained for simple domains of quadratic growth. 

We also note that the techniques used in the proof of Theorem \ref{thm: main} do not extend to higher GK dimension.  Nevertheless, we are unaware of a counter-example for algebras of GK dimension $3$.  We thus pose the following questions.
\begin{quest} Does there exist a non-noetherian finitely generated simple domain of GK dimension $3$?
\end{quest}
\begin{quest} Is the Krull dimension of a finitely generated simple noetherian algebra of GK dimension $3$ at most $2$?
\end{quest}
\bibliographystyle{amsplain}

\end{document}